\documentclass[12pt]{amsart}

\usepackage{latexsym}
\usepackage{psfrag}
\usepackage{amsmath}
\usepackage{amssymb}
\usepackage{epsfig}
\usepackage{amsfonts}

\newcommand{\excise}[1]{}

\newtheorem{thm}{Theorem}[section]

\newtheorem{conv}[thm]{Convention}

\newtheorem{Warn}[thm]{Caution}

    {\end{Example}}
    {\end{Remark}}

\def\sq{\square}

\def\rr{\mathbb R}

\def\ve{\varepsilon}

\def\ssu{\subset}

\def\<{\langle}
\def\>{\rangle}

\def\0{{\mathbf 0}}

\def\.{\hskip.06cm}
\def\ts{\hskip.03cm}

\def\vol{{\text {\rm vol}}}

\def\pt{\partial}

\begin{document}

\title[Inflating the cube]%
   {Inflating the cube without stretching}
\author[Igor~Pak]{Igor~Pak}
\date{25 July 2006}




\maketitle

\vspace{-1ex}

Consider the surface~$S$ of a unit cube~$C$ in~$\rr^3$.
Think of~$S$ being a container made of cardboard and full of water.
Now, can one add {\it more} water into container?
In other words, we are asking whether one can bend~$S$ without
tearing and stretching in such a way that the volume increases?
While one is tempted to say no, the following result gives a
positive answer:

\medskip

\noindent
{\bf Theorem.}  {\it There exists a non-convex polyhedron whose
surface is isometric to the surface of a cube of smaller volume.}

\medskip

Here by \emph{isometric} we mean that the geodesic distance between
pairs of points on the non-convex polyhedron is always equal to
the geodesic distance between of the corresponding pairs of points
on a cube.  Alternatively, it means that two surfaces can
be triangulated in such a way that they now consist of congruent
triangles which are glued according to the same combinatorial rules.
For example, if we push a vertex~$v$ of a cube~$C$ inside as shown
in the Figure below,
we obtain a polyhedron~$P$ whose surface~$\partial P$ is isometric
to~$S =\partial C$.  Of course,  $\vol(P) < \vol(C)$ in this case.

\begin{figure}[hbt]
\begin{center}
\psfrag{B}{$P$}
\psfrag{S}{$C$}
\psfrag{v}{$v$}
\epsfig{file=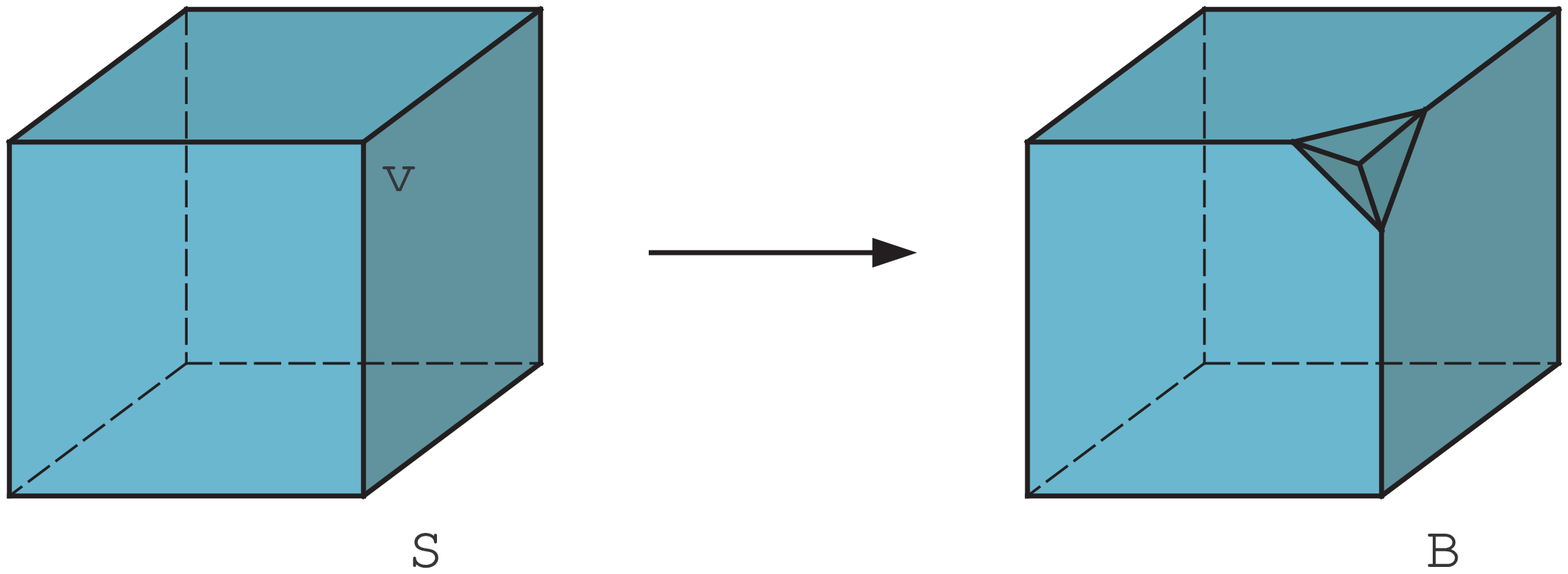,width=7.2cm}
\end{center}
\end{figure}

The theorem is somewhat unusual in a sense that one tends
to assume that convex objects always maximize the volume.
On the other hand, the Alexandrov uniqueness
theorem~\cite{A1} states that convex polyhedra
are uniquely determined (up to a rigid motion)
by the intrinsic geometry of their surfaces.
This means that \emph{all} polyhedra with surfaces
isometric to~$S$ are either congruent to~$C$ or
necessarily non-convex.


Before we present the proof of the theorem, let us say a few words
about the history of the problem.  Define a \emph{bending} to be
a continuous isometric deformation of a surface.
About ten years ago, two
different constructions of volume-increasing bending have
appeared~\cite{B,BZ} (see also~\cite{A2,S}).
In fact, Bleecker~\cite{B} showed that one can start with any convex
polytope~$P$ and bend it into a (non-convex) polyhedron of larger
volume.  Unfortunately, both constructions are technically involved.
Below we present a simple volume-increasing bending of a unit cube.
Our construction is based on the work of Milka~\cite{M}, where
all symmetric bendings of regular polyhedra were classified.

\medskip
\noindent
{\it Proof of the Theorem.} \.
Consider a cube~$C$ with side length~$1$ and the surface
$S = \partial C$.  Fix a parameter
$\ve \in (0,\frac12)$.  Think of~$\ve$ as being very small.
On each face of a cube, from every corner remove a square
of size-length~$\ve$.  Denote by~$R \ssu S$ the resulting
surface with boundary.
On every face of the cube there are four
boundary points which form a square.  Call these points corners.
Move each square away from the center of the cube (without
extending it or breaking the symmetry) until all of the distances
between the nearest corners on adjacent faces reach~$2\ts\ve$.
Take the convex hull of the corners  to obtain a polytope~$Q_\ve$.
To each triangular face of~$Q_\ve$ attach a triangular pyramid
whose base is equilateral with side-length~$2\ts\ve$
and whose other faces are right triangles.
Denote by~$P_\ve$ the resulting (non-convex) polyhedron.

Note that the surface~$\pt Q_\ve$ without triangular faces
between the corners is isometric to~$R$.  Similarly, three
$\ve \time \ve$ squares meeting at a vertex of the
cube can be bent into three faces of a pyramid
(see the Figure below). This easily implies that the
surface~$\pt P_\ve$ is isometric to the surface~$S$.

\begin{figure}[hbt]
\begin{center}
\psfrag{8}{{\bf 8}}
\epsfig{file=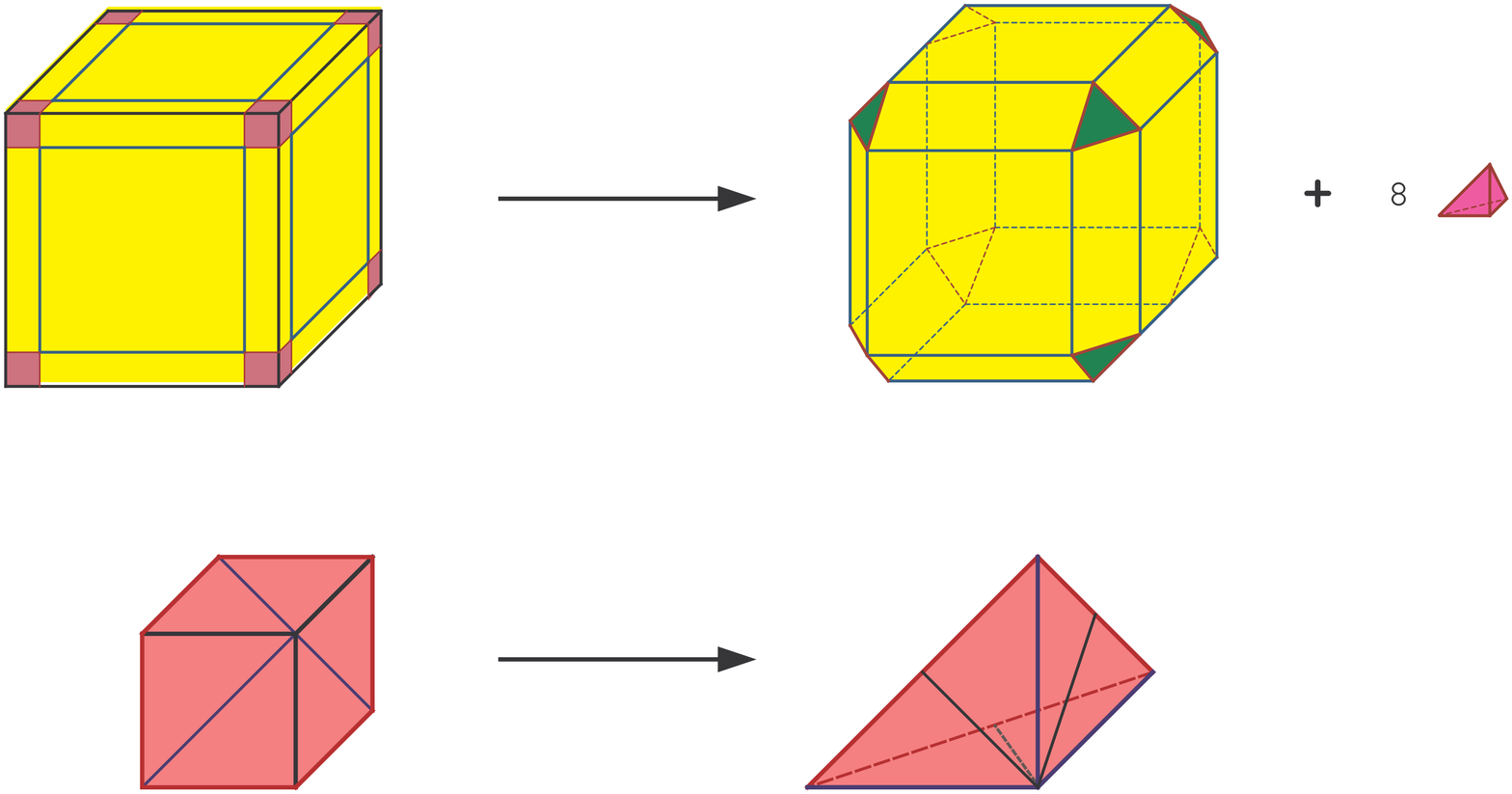,width=11.cm}
\end{center}
\end{figure}

Let us calculate the volume of the polyhedron~$P_\ve$.  Cut~$Q_\ve$
with six planes, each parallel to a square face and containing
the nearest edges of its four neighborhood square faces.
This subdivides~$Q_\ve$ into
one (interior) cube, six slabs (along the faces), twelve
right triangular prisms (along the edges),
and eight pyramids (one per cube vertex).
Observe that the cutting planes are at distance
$d=2\ts\ve/\sqrt{2} = \sqrt{2}\,\ve$ from the sides
of the interior cube.  We have:
$$
\vol(Q_\ve) \, = \, (1-2\ve)^3 + 6\cdot (1-2\ve)^2 d +
12\cdot(1-2\ve)\.\frac{d^2}{2} + 8\cdot \frac{d^3}{6}\,.
$$
Since $\vol(P_\ve) = \vol(Q_\ve) + 8 \cdot d^3/6$, we conclude
that
$$\vol(P_\ve) \, = \, 1 \, + \, 6\ts(\sqrt{2}-1)\. \ve \, + \, k_2 \. \ve^2
\, + \, k_3 \. \ve^3\.,
$$
where~$k_2$ and~$k_3$ are constants.
For small $\ve>0$, the terms $1+ 6(\sqrt{2}-1)\ve$ dominate, so
that $\vol(P_\ve)>1 = \vol(C)$, as desired. \ $\sq$

\bigskip

\noindent
{\bf Remark.} \.
That the polyhedra~$P_\ve$ are indeed non-convex may not be
immediately obvious.  As we mentioned earlier, this follows
from the Alexandrov theorem.  We leave the checking to the reader.

When $\ve \to \frac12$, the polyhedra~$Q_\ve$ converge to a regular
octahedron, and~$P_\ve$ converges to a nice stellated polyhedron
consisting of~18 right triangular faces.  The volume of the latter
is~$<0.95$, as the reader can deduce from the calculations above.
Thus at some point the volume of~$P_\ve$ stops increasing.

\bigskip

\noindent
{\bf Acknowledgements.} \. The author is grateful to David Jerison
and Ezra Miller for interesting discussions.  The author was
partially supported by the NSF.

\vskip.9cm


\vskip.6cm

\noindent
\textbf{\textsl{Igor Pak}}

\noindent
\textsl{Department of Mathematics}

\noindent
\textsl{M.I.T.}

\noindent
\textsl{Cambridge, MA 02139}

\noindent
\texttt{pak@math.mit.edu}

\end{document}